\documentclass{amsart}
\newtheorem{theorem}{Theorem}[section]
\newtheorem{lemma}[theorem]{Lemma}
\newtheorem{corollary}[theorem]{Corollary}

\theoremstyle{definition}

\theoremstyle{remark}
\newtheorem{remark}[theorem]{Remark}
\newtheorem{example}[theorem]{Example}

\newcommand{\cA}{{\mathcal A}}

\newcommand{\cD}{{\mathcal D}}

\newcommand{\cP}{{\mathcal P}}

\newcommand{\bC}{{\mathbb{C}}}
\newcommand{\bN}{{\mathbb{N}}}
\newcommand{\CP}{{\mathcal{CP}}}
\newcommand{\CcP}{{\mathcal C}c{\mathcal P}}
\numberwithin{equation}{section}

\newcommand{\jed}{{\mathbb{I}}}

\title[Positive maps]{Decomposability of extremal positive unital maps on $M_2(\bC)$}
\author[W. A. Majewski]{W{\l}adys{\l}aw A. Majewski}
\address{Institute of Theoretical Physics and Astrophysics, Gda{\'n}sk University,
Wita Stwosza 57, 80-952 Gda{\'n}sk, Poland}
\email{fizwam@univ.gda.pl}
\author[M. Marciniak]{Marcin Marciniak}
\address{Institute of Mathematics, Gda{\'n}sk University,
Wita Stwosza 57, 80-952 Gda{\'n}sk, Po\-land}
\email{matmm@univ.gda.pl}
\thanks{W.A.M. was partially supported by the grant
PBZ-MIN-008/PO3/2003, while M.M. was partially supported by
KBN grant 2P04A00723. The authors would like also thanks the
the support of EU RTN HPRN-CT-2002-00729 and
Poland-South Africa Cooperation Joint project.}
\begin{document}
\begin{abstract}
A map $\varphi:M_m(\bC)\to M_n(\bC)$ is decomposable if it is of the
form $\varphi=\varphi_1+\varphi_2$ where $\varphi_1$ is a CP map
while $\varphi_2$ is a co-CP map. It is known 
that if $m=n=2$ then every positive map is decomposable. Given an
extremal unital positive map $\varphi:M_2(\bC)\to M_2(\bC)$ we
construct concrete maps (not necessarily unital) $\varphi_1$ and $\varphi_2$ which give a
decomposition of $\varphi$. We also show that in most cases this
decomposition is unique.
\end{abstract}
\maketitle
\section{Introduction}
If $\cA$ is a $C^*$-algebra and $n\in\bN$ then by $M_n(\cA)$ we
denote the $C^*$-algebra of square $n\times n$-matrices with
coefficients in $\cA$. In particular $M_n(\bC)$ is the algebra of
matrices with complex entries. For each $m,n\in\bN$ we have the
following isomorphisms:
\begin{equation}\label{Isomorphism}
M_m(M_n(\bC))\cong M_m(\bC)\otimes M_n(\bC)\cong M_{mn}(\bC).
\end{equation}
It follows that $M_m(M_n(\bC))$ has the natural structure of a
$C^*$-algebra. In particular one defines the conjugation of
$\mathbf{A}=[A_{ij}]_{i,j=1}^m\in M_m(M_n(\bC))$ (where $A_{ij}\in
M_n(\bC)$ for $i,j=1,\ldots,m$) by the formula
$\mathbf{A}^*=[A_{ji}^*]_{i,j=1}^m$. Recall (see for example
\cite{T,ER}) that $\mathbf{A}$ is positive in $M_m(M_n(\bC))$ if and
only if $\sum_{i,j=1}^m\overline{\mu_i}\mu_j\langle
v_i,A_{ij}v_j\rangle\geq 0$ (i.e. the matrix $[\langle
v_i,A_{ij}v_j\rangle]_{i,j=1}^m$ is positive element of $M_m(\bC)$)
for every $v_1,\ldots,v_m\in\bC^n$ and $\mu_1,\ldots,\mu_m\in\bC$.
We say that $\mathbf{A}$ is \textit{block-positive} if
$\sum_{i,j=1}^m\overline{\mu_i}\mu_j\langle v,A_{ij}v\rangle\geq 0$
(i.e. the matrix $[\langle v,A_{ij}v\rangle]_{i,j=1}^m$ is positive
in $M_m(\bC)$) for every $v\in\bC^n$ and $\mu_1,\ldots,\mu_m\in\bC$.
For every $\mathbf{A}=[A_{ij}]_{i,j=1}^m\in M_m(M_n(\bC))$ we define
the \textit{partial transposition} of $\mathbf{A}$ by
$\mathbf{A}^\tau=[A_{ji}]_{i,j=1}^m$. Note the difference between
this operation and the usual transposition
$\mathbf{A}\mapsto\mathbf{A}^T$ on the algebra $M_{mn}(\bC)$
(\textit{cf.} (\ref{Isomorphism})): the usual transposition
preserves positivity of $\mathbf{A}$ while, for $n,m\geq 2$, partial
transposition does not!

A linear map $\varphi:M_m(\bC)\longrightarrow M_n(\bC)$ is called a
\textit{positive map} if $\varphi(A)$ is a positive matrix for every
positive matrix $A\in M_m(\bC)$. If $k\in\bN$ then $\varphi$ is
called \textit{$k$-positive map} (respectively
\textit{$k$-copositive map}) whenever $[\varphi(A_{ij})]_{i,j=1}^k$
(respectively $[\varphi(A_{ji})]_{i,j=1}^k$) is a positive element
in the algebra $M_k(M_n(\bC))$ for every positive element
$[A_{ij}]_{i,j=1}^k$ from $M_k(M_m(\bC))$. If $\varphi$ is
$k$-positive (respectively $k$-copositive) for every $k\in\bN$ then
$\varphi$ is called \textit{completely positive} or CP (respectively
\textit{completely copositive} or co-CP). A positive map which is a
sum of completely positive and completely copositive maps is called
\textit{decomposable}\footnote{In this paper we follow the
definition of decomposability given by St{\o}rmer in \cite{Sto}.
Note that there is no connection of this notion with decomposable
maps considered by U. Haagerup in the theory of operator spaces (see
\cite{ER}).}.

Let $\{E_{ij}\}_{i,j=1}^m$ be a system of matrix units in $M_m(\bC)$
and $\mathbf{H}_\varphi=[\varphi(E_{ij})]_{i,j=1}^m\in
M_m(M_n(\bC))$ be \textit{Choi matrix} of $\varphi$ with respect to
the system $\{E_{ij}\}$ (\cite{Choi}, see also \cite{MM}). Recall
the following
\begin{theorem}[\cite{Choi}, see also \cite{MM}]\label{Choimatrix}
Let $\varphi:M_m(\bC)\longrightarrow M_n(\bC)$ be a linear map. Then
\begin{enumerate}
\item the map $\varphi$ is positive if and only if the matrix
$\mathbf{H}_\varphi$ is block-positive;
\item the map $\varphi$ is completely positive
(respectively completely copositive) if and only if $\mathbf{H}_\varphi$
(respectively $\mathbf{H}_\varphi^\tau$) is
a positive element of $M_m(M_n(\bC))$.
\end{enumerate}
\end{theorem}

We say that a positive map $\varphi$ is \textit{unital} if
$\varphi(\jed)=\jed$ where $\jed$ denotes the identity matrix of the
respective algebra. The set of all positive (respectively completely
positive, completely copositive, decomposable) maps from $M_m(\bC)$
into $M_n(\bC)$ will be denoted by $\cP(m,n)$ (respectively
$\CP(m,n)$, $\CcP(m,n)$, $\cD(m,n)$).
We will write simply $\cP$, $\CP$, $\CcP$ and $\cD$ instead of
$\cP(m,n)$, $\CP(m,n)$, $\CcP(m,n)$ and $\cD(m,n)$ when no confusion
will arise. Observe that all of these sets have the structure of a
convex cone. By $\cP_1$, $\CP_1$, $\CcP_1$ and $\cD_1$ we will
denote the subsets of unital maps from respective cones. All of them
are convex subsets.

Let $C$ be a convex cone and $c\in C$. We say that $c$ is an extreme
point of $C$ if for every $c_1,c_2\in C$ the equality $c=c_1+c_2$
implies $c_1=\lambda c$ and $c_2=(1-\lambda)c$ for some
$0\leq\lambda\leq 1$. The generalization of the notion of
extremality leads to the face structure of the cone $C$. Namely, we
say that a convex subcone $F\subset C$ is a \textit{face} of $C$ if
for every $c_1,c_2\in C$ the condition $c_1+c_2\in F$ implies
$c_1,c_2\in F$. Kye in \cite{Kye} gave an interesting
characterization of maximal faces of the cone $\cP(m,n)$
\begin{theorem}[\cite{Kye}]\label{maxface}
A convex subset $F\subset\cP(m,n)$ is a maximal face of $\cP(m,n)$
if and only if there are vectors $\xi\in\bC^m$ and $\eta\in\bC^n$
such that $F=F_{\xi,\eta}$ where
\begin{equation}\label{face}
F_{\xi,\eta}=\{\varphi\in\cP(m,n):\,\varphi(P_\xi)\eta=0\}
\end{equation}
and $P_\xi$ denotes the one-dimensional orthogonal projection in
$M_m(\bC)$ onto the subspace generated by the vector $\xi$.
\end{theorem}

\section{The case $m=n=2$}
In this section we analyze in details the case $m=n=2$. In
\cite{Sto} the following characterization of extremal point of
$\cP_1$ is given
\begin{theorem}\label{Stormer}
A positive unital map $\varphi:M_2(\bC)\longrightarrow M_2(\bC)$ is
an extremal point of $\cP_1$ if and only if there are unitary
operators $V,W\in U(2)$ such that the Choi matrix of the map
$\varphi_{V,W}:A\mapsto V^*\varphi(WAW^*)V$ has the form
\begin{equation}\label{extrem}
\left[\begin{array}{cc|cc}
1&0&0&y\\0&b&\overline{z}&t\\ \hline 0&z&0&0\\ \overline{y}&\overline{t}&0&u
\end{array}\right]
\end{equation}
where the coefficients fulfil the following relations:
\begin{enumerate}
\item $b\geq 0$, $u\geq 0$ and $b+u=1$,
\item $|t|^2=2b(u-|y|^2-|z|^2)$ in the case when $b\neq 0$,
and $|y|=1$ or $|z|=1$ when $b=0$.
\end{enumerate}
\end{theorem}

By $\{e_1,e_2\}$ we denote the canonical basis in $\bC^2$. Let
$\varphi$ be an extremal positive unital map as in Theorem
\ref{Stormer}. One can observe that $\varphi\in F_{\xi,\eta}$ for
$\xi=We_2$ and $\eta=Ve_1$ where $V,W$ are the unitary operators
from Theorem \ref{Stormer}. Suppose that
$\varphi=\varphi_1+\varphi_2$ where $\varphi_1$ is a completely
positive map while $\varphi_2$ is a completely copositive one. Then
both $\varphi_1$ and $\varphi_2$ should be elements of
$F_{\xi,\eta}$ because $F_{\xi,\eta}$ is a face (\textit{cf.}
Theorem \ref{maxface}).
\begin{remark}
There are extremal maps of the form (\ref{extrem}) which are neither
completely positive nor completely copositive (see Example
\ref{examp} below). On the other hand theorem of Woronowicz (see
\cite{Wor}) asserts that every map from $\cP(2,2)$ is decomposable.
Hence, the maps $\varphi_1$ and $\varphi_2$ giving the decomposition
of an extremal element of $\cP_1(2,2)$ do not need be scalar
multiples of unital maps.
\end{remark}

In the sequel we will use the following lemmas.
\begin{lemma}\label{Choifacelemma}
Let $\psi\in F_{\xi,\eta}$ for some $\xi,\eta\in\bC^2$, and $V$ and
$W$ be unitary operators from $M_2(\bC)$ such that $\xi=We_2$ and
$\eta=Ve_1$. Then the Choi matrix of the map $\psi_{V,W}$ has the
form
\begin{equation}\label{Choiface}
\left[\begin{array}{cc|cc}
a&c&0&y\\\overline{c}&b&\overline{z}&t\\\hline 0&z&0&0\\\overline{y}&\overline{t}&0&u
\end{array}\right]
\end{equation}
for some $a,b,u\geq 0$ and $c,y,z,t\in\bC$. Moreover, the following conditions hold:
\begin{enumerate}
\item $|c|^2\leq ab$,
\item $|t|^2\leq bu$,
\item $(|y|+|z|)^2\leq au$.
\end{enumerate}
\end{lemma}
\begin{proof}
Let us write shortly $\psi'$ instead of $\psi_{V,W}$. The Choi
matrix of the map $\psi'$ has the form
$$\mathbf{H}=\left[\begin{array}{cc}
\psi'(E_{11}) & \psi'(E_{12}) \\ \psi'(E_{21}) & \psi'(E_{22})
\end{array}\right].$$
From the condition (\ref{face}) we get
$$\psi'(E_{22})e_1=V^*\psi(WE_{22}W^*)Ve_1=V^*\psi(P_\xi)\eta=0.$$
Because $\psi'(E_{22})$ is a positive element of $M_2(\bC)$ then
$\psi'(E_{22})=uE_{22}$ for some $u\geq 0$. Hence
$\psi'(E_{22})=\left[\begin{array}{cc}0&0\\0&u\end{array}\right]$.
Positivity of $\psi'$ implies also that $\psi'(E_{11})$ is a
hermitian matrix of the general form
$\psi'(E_{11})=\left[\begin{array}{cc}a&c\\\overline{c}&b\end{array}\right]$
with $a,b\geq 0$ and
$\det\psi'(E_{11})=ab-|c|^2\geq 0$, so we have (1). Let $\psi'(E_{12})=\left[\begin{array}{cc}x&y\\
\overline{z}&t\end{array}\right]$
for some $x,y,z,t\in\bC$. Because $\psi'$ is positive then
$\psi'(E_{21})=\psi'(E_{12}^*)=\psi'(E_{12})^*=\left[\begin{array}{cc}\overline{x}&z\\
\overline{y}&\overline{t}\end{array}\right]$. By Theorem
\ref{Choimatrix} $\mathbf{H}$ is block-positive, hence the matrix
$$\left[\begin{array}{cc}\langle e_1,\psi'(E_{11})e_1\rangle &
\langle e_1,\psi'(E_{12})e_1\rangle \\
\langle e_1,\psi'(E_{21})e_1\rangle & \langle
e_1,\psi'(E_{22})e_1\rangle \end{array}\right]=
\left[\begin{array}{cc} a&x\\\overline{x}&0\end{array}\right]$$ is
positive and consequently $x=0$. Thus we arrived to the form
(\ref{Choiface}). Another application of block-positivity of $H$
leads to the conclusion that the matrix
$$\left[\begin{array}{cc}\langle e_2,\psi'(E_{11})e_2\rangle & \langle e_2,\psi'(E_{12})e_2\rangle \\
\langle e_2,\psi'(E_{21})e_2\rangle & \langle e_2,\psi'(E_{22})e_2\rangle \end{array}\right]=
\left[\begin{array}{cc} b&t\\\overline{t}&u\end{array}\right]$$
is positive, hence we get inequality (2).

Let $\omega$ be a linear functional on $M_2(\bC)$. By Corollary 8.4
in \cite{Sto} $\omega$ is a positive functional if and only if
$\omega(E_{11})\geq 0$, $\omega(E_{22})\geq 0$,
$\omega(E_{21})=\overline{\omega(E_{12})}$ and
$|\omega(E_{12})|^2\leq\omega(E_{11})\omega(E_{22})$. Let us denote
$\alpha=\omega(E_{11})$, $\beta=\omega(E_{22})$ and
$\gamma=\omega(E_{12})$. Obviously $\omega\circ\psi'$ is a positive
functional for every positive $\omega$. From another application of
this result of St{\o}rmer we get that the inequality
\begin{equation}\label{posit1}
|\gamma y+\overline{\gamma}\overline{z}+\beta t|^2\leq
\beta u(\alpha a+\beta b+2\Re(\gamma c))
\end{equation}
holds whenever $\alpha\geq 0$, $\beta\geq 0$ and
$|\gamma|^2\leq\alpha\beta$. Inequality (\ref{posit1}) can be
written in the form
\begin{equation}\label{posit2}
|\gamma y+\overline{\gamma}\overline{z}|^2+\beta^2|t|^2+
2\Re[(\gamma y+\overline{\gamma}\overline{z})\beta
\overline{t}]\leq
\beta u(\alpha a+\beta b+2\Re(\gamma c)).
\end{equation}
Putting here $-\gamma$ instead of $\gamma$ we get
\begin{equation}\label{posit3}
|\gamma y+\overline{\gamma}\overline{z}|^2+\beta^2|t|^2-
2\Re[(\gamma y+\overline{\gamma}\overline{z})
\beta \overline{t}]\leq
\beta u(\alpha a+\beta b-2\Re(\gamma c)).
\end{equation}
If we add (\ref{posit2}) and (\ref{posit3}) and divide the result by
$2$, then we get
\begin{equation}\label{posit4}
|\gamma y+\overline{\gamma}\overline{z}|^2+\beta^2|t|^2\leq
\beta u(\alpha a+\beta b).
\end{equation}
This can be rewritten equivalently as
\begin{equation}\label{posit5}
|\gamma|^2(|y|^2+|z|^2)+2\Re(yz\gamma^2)+\beta^2|t|^2\leq
\beta u(\alpha a+\beta b)
\end{equation}
Let $\varepsilon>0$ and take $\alpha=\varepsilon^{-1}$,
$\beta=\varepsilon$ and $\gamma$ such that $|\gamma|=1$ and
$yz\gamma^2=|y||z|$. Then (\ref{posit5}) has the form
\begin{equation}\label{posit6}
(|y|+|z|)^2\leq au +\varepsilon^2(bu-|t|^2)
\end{equation}
Because $\varepsilon$ can be chosen arbitrary small then we get the
inequality (3) and the proof is finished.
\end{proof}
\begin{lemma}\label{ChoiCP}
A map $\psi\in F_{\xi,\eta}$ is completely positive if and only if
the coefficients of the matrix $\mathbf{H}$ from (\ref{Choiface})
fulfil the following conditions:
\begin{enumerate}
\item[(A1)] $z=0$,
\item[(A2)] $|y|^2\leq au$,
\item[(A3)] $|t|^2\leq bu$,
\item[(A4)] $|c|^2\leq ab$,
\item[(A5)] $a|t|^2+u|c|^2\leq b(au-|y|^2)+2\Re (ct\overline{y})$.
\end{enumerate}
Analogously, $\psi$ is completely copositive if and only if the
following conditions hold:
\begin{enumerate}
\item[(B1)] $y=0$,
\item[(B2)] $|z|^2\leq au$,
\item[(B3)] $|t|^2\leq bu$,
\item[(B4)] $|c|^2\leq ab$,
\item[(B5)] $a|t|^2+u|c|^2\leq b(au-|z|^2)+2\Re (c\overline{t}\overline{z})$.
\end{enumerate}
\end{lemma}
\begin{proof}
By Theorem \ref{Choimatrix} and properties of unitary equivalence
$\psi$ is completely positive if and only if the matrix $\mathbf{H}$
is positive. This is equivalent to the fact that all principal
minors of $\mathbf{H}$ are nonnegative. Conditions (A2), (A3) and
(A4) follow from the fact that
$\left|\begin{array}{cc}a&y\\\overline{y}&u\end{array}\right|\geq
0$,
$\left|\begin{array}{cc}b&t\\\overline{t}&u\end{array}\right|\geq 0$
and
$\left|\begin{array}{cc}a&c\\\overline{c}&b\end{array}\right|\geq
0$. (A1) is a consequence of the equality
$\det\mathbf{H}=-|z|^2(au-|y|^2)$ and (A2). Inequality in (A5) is
equivalent to
$\left|\begin{array}{ccc}a&c&y\\\overline{c}&b&t\\\overline{y}&\overline{t}&u\end{array}\right|\geq
0$.

The second part of the lemma follows in the similar way from
positivity of the matrix $\mathbf{H}^\tau$ in the case when $\psi'$
is completely copositive.
\end{proof}
\begin{remark}
If $\varphi$ is an extremal positive unital map described in Theorem
\ref{Stormer} with the Choi matrix (\ref{extrem}) then by condition
(3) from Lemma \ref{Choifacelemma} $|y|+|z|\leq u^{1/2}$. Lemma 8.11
in \cite{Sto} claims that in the case $b>0$ a stronger condition
holds. Namely,
\begin{equation}\label{yz}
|y|+|z|=u^{1/2}.
\end{equation}
Moreover, it follows from Lemma 8.8 in \cite{Sto} that in this case
\begin{equation}\label{t}
t^2=-4(1-u)y\overline{z}
\end{equation}
\end{remark}
\begin{example}\label{examp}
Consider the map $\psi:M_2(\bC)\longrightarrow M_2(\bC)$ with the
Choi matrix
$$\left[\begin{array}{cc|cc}
1&0&0&\frac{1}{2}s\\0&1-s^2&\frac{1}{2}s&is(1-s^2)^{1/2}\\\hline
0&\frac{1}{2}s&0&0\\\frac{1}{2}s&-is(1-s^2)^{1/2}&0&s^2\end{array}\right]$$
where $0<s<1$. It follows from Theorem \ref{Stormer} (compare also
with (\ref{yz}) and (\ref{t})) and Lemma \ref{ChoiCP} that $\psi$ is
an extremal positive unital map which is neither completely positive
nor completely copositive.
\end{example}

Now, we are ready to formulate our main theorem
\begin{theorem}\label{main}
Assume that $\varphi$ is an extremal positive unital map with the
Choi matrix given in (\ref{extrem}) and $u>0$, $y\neq 0$, $z\neq 0$.
Then there are $\varphi_1,\varphi_2\in F_{\xi,\eta}$ such that
$\varphi_1$ is completely positive, $\varphi_2$ is completely
copositive and $\varphi=\varphi_1+\varphi_2$. Moreover, the pair
$\varphi_1,\varphi_2$ is uniquely determined, and Choi matrices
$\mathbf{H}_1$, $\mathbf{H}_2$ of maps $\varphi_1$, $\varphi_2$ have
the following form:
\begin{equation}\label{CP}
\mathbf{H}_1=\left[\begin{array}{cc|cc} |y|u^{-1/2} &
c 
& 0 & y
\\
\overline{c}
& |z|(1-u)u^{-1/2} & 0 & \frac{1}{2}t \\\hline 0&0&0&0
\\ \overline{y} & \frac{1}{2}\overline{t} & 0 &
|y|u^{1/2}\end{array}\right],
\end{equation}
\begin{equation}\label{CcP}
\mathbf{H}_2=\left[\begin{array}{cc|cc} |z|u^{-1/2} &
-c
& 0 & 0
\\
-\overline{c}
& |y|(1-u)u^{-1/2} & \overline{z} & \frac{1}{2}t
\\\hline 0&z&0&0 \\ 0 &
\frac{1}{2}\overline{t} & 0 &
|z|u^{1/2}\end{array}\right]
\end{equation}
where $u,y,z,t$ are the coefficients of the matrix (\ref{extrem})
and $c$ is a complex number such that $c^2=-(1-u)u^{-1}yz$.
\end{theorem}
\begin{remark}
The uniqueness of the decomposition of $\varphi$ onto $\varphi_1$
and $\varphi_2$ does not hold if the assumptions of the above
theorem are not fulfilled. To see this let us consider the following
cases:
\begin{enumerate}
\item
$u=0$. Then the Choi matrix of the map $\varphi$ has the following
form $$\mathbf{H}=\left[\begin{array}{cc|cc}1&0&0&0\\0&1&0&0\\\hline
0&0&0&0\\0&0&0&0\end{array}\right]$$ and $\varphi$ is completely
positive and completely copositive.
\item
$y=0$. The Choi matrix of $\varphi$ has the form
$$\mathbf{H}=\left[\begin{array}{cc|cc}1&0&0&0\\0&1-|z|^2&\overline{z}&0\\\hline
0&z&0&0\\0&0&0&|z|^2\end{array}\right]$$ and the map $\varphi$ is
completely copositive.
\item
$z=0$. The Choi matrix of $\varphi$ has the form
$$\mathbf{H}=\left[\begin{array}{cc|cc}1&0&0&y\\0&1-|y|^2&0&0\\\hline
0&0&0&0\\\overline{y}&0&0&|y|^2\end{array}\right]$$ and the map
$\varphi$ is completely positive.
\end{enumerate}
Turning to non-uniqueness of the decomposition of $\varphi$, let
$\varepsilon>0$ and put
$$\mathbf{A}_\varepsilon=\left[\begin{array}{cc|cc}0&0&0&0\\
0&\varepsilon&0&0\\\hline 0&0&0&0\\0&0&0&0\end{array}\right].$$
Then, it follows from Lemma \ref{ChoiCP} that
$\mathbf{A}_\varepsilon$ determines the map $\psi_\varepsilon$ which
is completely positive and completely copositive. Moreover, in each
of the above three cases the equality
$\mathbf{H}=\left(\mathbf{H}-\mathbf{A}_\varepsilon\right)+\mathbf{A}_\varepsilon$
describes the decomposition onto CP and co-CP parts for every
sufficiently small $\varepsilon$.
\end{remark}
\begin{proof}[Proof of Theorem \ref{main}]
The existence of the decomposition follows from \cite{MM,Sto} and
the discussion after Theorem \ref{Stormer}. To show that the
decomposition is unique (and has the required form) assume that it
is given by $\mathbf{H}=\mathbf{H}_1+ \mathbf{H}_2$, where
\begin{equation}\label{sum}
\mathbf{H}_1=\left[\begin{array}{cc|cc}a_1&c&0&y\\\overline{c}&b_1&0&t_1\\
\hline 0&0&0&0\\\overline{y}&\overline{t_1}&0&u_1
\end{array}\right],\;\;\;\mathbf{H}_2=\left[\begin{array}{cc|cc}
a_2&-c&0&0\\-\overline{c}&b_2&\overline{z}&t_2\\\hline 0&z&0&0\\0&\overline{t_2}&0&u_2
\end{array}\right].
\end{equation}
Then the above coefficients fulfil the following set of relations:
\begin{equation}\label{suma}
a_1+a_2=1,
\end{equation}
\begin{equation}\label{sumb}
b_1+b_2=1-u,
\end{equation}
\begin{equation}\label{sumt}
t_1+t_2=t,
\end{equation}
\begin{equation}\label{sumu}
u_1+u_2=u,
\end{equation}
\begin{equation}\label{CP1}
|y|^2\leq a_1u_1,
\end{equation}
\begin{equation}\label{CP2}
|t_1|^2\leq b_1u_1,
\end{equation}
\begin{equation}\label{CP3}
|c|^2\leq a_1b_1,
\end{equation}
\begin{equation}\label{CP4}
a_1|t_1|^2+u_1|c|^2\leq b_1(a_1u_1-|y|^2)+2\Re(ct_1\overline{y}),
\end{equation}
\begin{equation}\label{CcP1}
|z|^2\leq a_2u_2,
\end{equation}
\begin{equation}\label{CcP2}
|t_2|^2\leq b_2u_2,
\end{equation}
\begin{equation}\label{CcP3}
|c|^2\leq a_2b_2,
\end{equation}
\begin{equation}\label{CcP4}
a_2|t_2|^2+u_2|c|^2\leq b_2(a_2u_2-|z|^2)-2\Re(c\overline{t_2}\overline{z}).
\end{equation}
We divide the rest of the proof into some lemmas.
\begin{lemma}\label{au}
Coefficients $a_1,a_2,u_1,u_2$ have the following form
\begin{eqnarray}
a_1&=&|y|u^{-1/2},\label{a1}\\
a_2&=&|z|u^{-1/2},\label{a2}\\
u_1&=&|y|u^{1/2},\label{u1}\\
u_2&=&|z|u^{1/2}.\label{u2}
\end{eqnarray}
\end{lemma}
\begin{proof}
From (\ref{yz}) we have $|y|u^{-1/2}+|z|u^{-1/2}=1$. Let
$p=|y|u^{-1/2}$, $q=a_1$ and $r=u_1u^{-1}$. Then (\ref{CP1}) gives
\begin{equation}\label{sys1}
p^2\leq q r
\end{equation}
while (\ref{yz}), (\ref{CcP1}), (\ref{suma}) and (\ref{sumu}) lead
to
\begin{equation}\label{sys2}
(1-p)^2\leq (1-q)(1-r).
\end{equation}
The system of inequalities (\ref{sys1}) and (\ref{sys2}) is
equivalent to
\begin{equation}\label{sys}
\frac{p^2}{q}\leq r\leq 1-\frac{(1-p)^2}{1-q}.
\end{equation}
By simple calculations one can show that the inequality between the
first and the last terms in (\ref{sys}) is equivalent to
$(q-p)^2\leq 0$, so $q=p$. So, puting $p$ instead of $q$ in
(\ref{sys}) we obtain also $r=p$. Hence, we have
\begin{eqnarray*}
&&a_1=q=p=|y|u^{-1/2},\\
&&a_2=1-a_1=1-|y|u^{-1/2}=|z|u^{-1/2},\\
&&u_1=r u=pu=|y|u^{1/2},\\
&&u_2=u-u_1=(u^{1/2}-|y|)u^{1/2}=|z|u^{1/2}.
\end{eqnarray*}
\end{proof}
\begin{lemma}\label{ct}
The following relations hold
\begin{equation}\label{ct1}
|y|t_1=y\overline{c}u^{1/2},
\end{equation}
\begin{equation}\label{ct2}
|z|t_2=-\overline{z}cu^{1/2}.
\end{equation}
\end{lemma}
\begin{proof}
Observe that application of (\ref{a1})-(\ref{u2}) reduces
inequalities (\ref{CP4}) and (\ref{CcP4}) to
\begin{equation}\label{f1}
|y|u^{-1/2}|t_1|^2+|y|u^{1/2}|c|^2-2\Re(ct_1\overline{y})\leq 0
\end{equation}
and
\begin{equation}\label{f2}
|z|u^{-1/2}|t_2|^2+|z|u^{1/2}|c|^2+2\Re(c\overline{t_2}\overline{z})\leq 0
\end{equation}
respectively. Let $y_1,z_1\in\bC$ be such that $y_1^2=y$ and
$z_1^2=z$. Then (\ref{f1}) and (\ref{f2}) can be rewritten in the
form
\begin{equation*}
\left\vert\overline{y_1}u^{-1/4}t_1-y_1u^{1/4}\overline{c}\right\vert^2\leq 0
\end{equation*}
and
\begin{equation*}
\left\vert\overline{z_1}u^{-1/4}\overline{t_2}+z_1u^{1/4}\overline{c}\right\vert^2\leq 0.
\end{equation*}
These inequalities are equivalent to
$\overline{y_1}u^{-1/4}t_1=y_1u^{1/4}\overline{c}$ and
$\overline{z_1}u^{-1/4}\overline{t_2}=-z_1u^{1/4}\overline{c}$.
Multiplication of both sides of the first equality by $y_1u^{1/4}$
leads to (\ref{ct1}) while multiplication of the second one by
$z_1u^{1/4}$ gives (\ref{ct2}).
\end{proof}
\begin{corollary}\label{cor}
$|t_1|=|t_2|\geq \frac{1}{2}|t|$.
\end{corollary}
\begin{proof}
The equality follows from (\ref{ct1}) and (\ref{ct2}) while the
inequality is a consequence of (\ref{sumt}) and the triangle
inequality.
\end{proof}
\begin{lemma}\label{b}
The following relations hold:
\begin{equation}\label{b1}
b_1=|z|(1-u)u^{-1/2},
\end{equation}
\begin{equation}\label{b2}
b_2=|y|(1-u)u^{-1/2}.
\end{equation}
\end{lemma}
\begin{proof}
From (\ref{t}), (\ref{CP2}), (\ref{CcP2}), (\ref{u1}), (\ref{u2})
and Corollary \ref{cor} we have the following inequalities
$$(1-u)|y||z|=\frac{1}{4}|t|^2\leq |t_1|^2\leq b_1u_1=|y|u^{1/2}b_1$$
and
$$(1-u)|y||z|=\frac{1}{4}|t|^2\leq |t_2|^2\leq b_2u_2=|z|u^{1/2}b_2.$$
From the first inequality we obtain
$$b_1\geq |z|(1-u)u^{-1/2}$$
while from the second one and (\ref{sumb}) we have
\begin{eqnarray*}
b_1&=&1-u-b_2\leq 1-u-|y|(1-u)u^{-1/2}=\\
&=& (u^{1/2}-|y|)(1-u)u^{-1/2}=|z|(1-u)u^{-1/2}.
\end{eqnarray*}
Thus we obtain (\ref{b1}). In a similar way we get (\ref{b2}).
\end{proof}
\begin{lemma}\label{tt}
$|t_1|=|t_2|=\frac{1}{2}|t|$.
\end{lemma}
\begin{proof}
It follows from (\ref{CP2}), Lemma \ref{b}, (\ref{u1}) and (\ref{t}) that
$$|t_1|^2\leq b_1u_1= |z|(1-u)u^{-1/2}\cdot |y|u^{1/2}=(1-u)|y||z|=\frac{1}{4}|t|^2.$$
The converse inequality is included in Corollary \ref{cor}.
\end{proof}
\begin{corollary}\label{cort}
$t_1=t_2=\frac{1}{2}t$.
\end{corollary}
\begin{proof}
It easily follows from (\ref{sumt}) and Lemma \ref{tt}.
\end{proof}
\begin{lemma}\label{c}
$c^2=-(1-u)u^{-1}yz$.
\end{lemma}
\begin{proof}
From (\ref{ct2}) and Corollary \ref{cort} we obtain
$c=-\frac{1}{2}z|z|^{-1}tu^{-1/2}$.
Thus, (\ref{t}) implies
$c^2=-z^2|z|^{-2}\cdot (1-u)y\overline{z}u^{-1}=-(1-u)u^{-1}yz$.
\end{proof}

The coefficient $c$ in (\ref{CP}) and (\ref{CcP}) is uniquely
determined. It can be described in the following way. Let
$y_1,z_1\in\bC$ be such that $y_1^2=y$, $z_1^2=z$ and
$t=2i(1-u)^{1/2}y_1\overline{z_1}$ (\textit{cf.} (\ref{t})). The
numbers $y_1$, $z_1$ are not uniquely determined by these conditions
but the expression $y_1z_1$ is. Then, by (\ref{ct2})
\begin{equation}\label{cex}
c=-i(1-u)^{1/2}u^{-1/2}y_1z_1.
\end{equation}

Summing results of Lemmas \ref{au}, \ref{b}, \ref{c} and Corollary
\ref{t} we end the proof of Theorem \ref{main}.
\end{proof}

\begin{corollary}
If $\varphi$ is extremal positive unital map with the Choi matrix of
the form (\ref{extrem}) and assumptions of Theorem \ref{main} are
fulfilled then $\varphi(A)=U_1AU_1^*+U_2A^TU_2^*$ for every $A\in
M_2(\bC)$, where $U_1,U_2\in M_2(\bC)$ are of the form
$$U_1=\left[\begin{array}{cc}y_1u^{-1/4}&0\\i\overline{z_1}(1-u)^{1/2}
u^{-1/4}&\overline{y_1}u^{1/4}\end{array}\right],\;\;\;
U_2=\left[\begin{array}{cc}z_1u^{-1/4}&0\\-i\overline{y_1}(1-u)^{1/2}
u^{-1/4}&\overline{z_1}u^{1/4}\end{array}\right],$$
where $y_1$ and $z_1$ are as in the proof of Theorem \ref{main}.
\end{corollary}

\vspace{3mm} \textit{Acknowledgments.} Part of the work was done
during a visit of the authors at University of South Africa in
Pretoria. The authors want to thank Louis E. Labuschagne for his
kind hospitality and fruitful discussions.

\end{document}